\documentclass[11pt]{article}
\usepackage{amscd}
\usepackage{amsmath}
\usepackage{latexsym}
\usepackage{amsfonts}
\usepackage{amssymb}
\usepackage{amsthm}
\usepackage{graphicx}

 \newtheorem{proposition}{Proposition}[section]
 
 \newtheorem{lemma}{Lemma}[section]
 \newtheorem{remark}{Remark}[section]
 \newtheorem{theorem}{Theorem}[section]

\begin{document}

 \title{\large \bf LOW REGULARITY GLOBAL WELL-POSEDNESS FOR THE KLEIN-GORDON-SCHR\"{O}DINGER
 SYSTEM WITH THE HIGHER ORDER YUKAWA COUPLING
 \footnote{This project is supported by the National Natural Science Foundation of
 China.}
      }
 \author{{Changxing Miao\ \ }\\
          {\small Institute of Applied Physics and Computational Mathematics}\\
         {\small P. O. Box 8009,\ Beijing,\ China,\ 100088}\\
         {\small (miao\_changxing@mail.iapcm.ac.cn) }\\ \\
         {Guixiang Xu\ \ } \\
        {\small The Graduate School of China Academy of Engineering Physics  }\\
        {\small P. O. Box 2101,\ Beijing,\ China,\ 100088 } \\
        {\small (gxxu78@hotmail.com) } \\ \\
         \date{}
        }
\maketitle

\begin{abstract}
In this paper, we consider the Klein-Gordon-Schr\"{o}dinger system
with the higher order Yukawa coupling in $ \mathbb{R}^{1+1} $, and
prove the local and global wellposedness in $L^2\times H^{1/2}$. The
method to be used is adapted from the scheme originally by
Colliander J., Holmer J., Tzirakis N. \cite{CoHT06} to use the
available $L^2$ conservation law of $u$ and control the growth of
$n$ via the estimates in the local theory.
\end{abstract}

 \begin{center}
 \begin{minipage}{120mm}
 \end{minipage}
\end{center}

\section{Introduction}
\setcounter{section}{1}\setcounter{equation}{0} The Cauchy problem
\begin{equation}\label{equ0}
\left\{ \aligned
   i u_t +\Delta u  = & - n u , \quad x \in \mathbb{R}^d, t\in \mathbb{R}; \\
   n_{tt} + (1-\Delta ) n  = & |u|^2 , \quad \ x \in \mathbb{R}^d, t\in \mathbb{R}; \\
   u(0)=u_0 ,\ \ n(0) = & n_0,\ \ n_t(0)=n_1.
\endaligned
\right.
\end{equation}
have been considered in \cite{BaC78}, \cite{FuT75a}, \cite{FuT75b},
\cite{FuT78}. Here $u : \mathbb{R}^d \times \mathbb{R}
\longrightarrow \mathbb{C}$ is the nucleon field and $n :
\mathbb{R}^d \times \mathbb{R} \longrightarrow \mathbb{R}$ is the
meson field. H. Pecher \cite{Pe04} considered the system
$(\ref{equ0})$ in $\mathbb{R}^{3+1}$ by Fourier truncation method
\cite{Bo98a}. N. Tzirakis \cite{Tz05b} consider the same system in
one, two, three dimension by I-method \cite{KeT98b}. They also
obtained a polynomial in time bound for the growth of the norms.
Recently, using the available $L^2$ conservation law of $u$ and
controlling the growth of $n$ via the estimate in the local theory,
J. Colliander et al. \cite{CoHT06} obtained the optimal global
well-posedness of $(\ref{equ0})$ in $\mathbb{R}^{3+1}$. It is also
applicable to 1D and 2D case.

Just as in \cite{Bi89}, \cite{Co98}, the system $(\ref{equ0})$ is
naturally generalized to the following system
\begin{equation*}
\left\{ \aligned
   i u_t +\Delta u  = & - n f(|u|^2)u ;\\
   n_{tt} + (1-\Delta ) n  = & F(|u|^2),\quad F'=f, F(0)=f(0)=0; \\
   u(0)=u_0 ,\ \ n(0) = & n_0 ,\ \ n_t(0)=n_1 .
\endaligned
\right.
\end{equation*}
The restricted case $F(s)=s^m$ in GKLS will be called KLS$_{m}$. In
this paper, we only consider the case $1 \leq m <2, d=1$, that is
\begin{equation}\label{equ1}
\left\{ \aligned
   i u_t +\partial^2_x u  = & - m n |u|^{2(m -1)}u , \quad x \in \mathbb{R}, t\in \mathbb{R};\\
   n_{tt} + (1-\partial^2_x ) n  = & |u|^{2m} , \quad  \qquad\qquad \quad x \in \mathbb{R}, t\in \mathbb{R}; \\
   (u, n, \dot{n})(0)=& (u_0, n_0, n_1) .
\endaligned
\right.
\end{equation}

The reason that the higher order powers are introduced into the
physically relevant dispersive PDEs is to adjust the strength of the
nonlinearity relative to the dispersion to work toward understanding
the balance between the two effects. We give the similar scaling
analysises in next section.

It is well known that the following conservation laws hold for
$(\ref{equ1})$:
\begin{equation}{\label{cl}}
\left\{ \aligned
M(u)(t) &=: \| u (t) \| ,\\
E(u, n)(t) &=: \|  \partial_x u (t) \|^2 + \displaystyle \frac{1}{2}
\big(\|A n (t)\|^2 +\|n_t(t)\|^2\big)
 -\displaystyle\int_{\mathbb{R}}
|u(t)|^{2m} n(t) dx
,\\
\endaligned
\right.
\end{equation}
where $\|\cdot\|$ denotes the norm of $L^2(\mathbb{R})$ and $A $
denotes $(I-\partial^2_{x})^{\frac{1}{2}}$. Here, we use the method
in \cite{CoHT06} rather than Fourier truncation method and I-method
to consider the low regularity. The idea is to use the available
$L^2$ conservation law of $u$ and control the growth of $n$ via the
estimates in the local theory.


Our main result is the following theorem
\begin{theorem}\label{the1}
Let $1 \leq m<2$, then the KLS$_{m}$ $(\ref{equ1})$ in dimension
$d=1$ is global well-posedness for $(u_0, n_0, n_1) \in L^2 \times
H^{1/2} \times H^{-1/2}$. More precisely, the solution $(u, n)\in
C(\mathbb R; L^2)\times C(\mathbb R; H^{\frac12})$ satisfies for
$t\in\mathbb R$,
$$\big\| u(t)\big\|_{L^2} = \big\| u_0\big\|_{L^2} $$
and
\begin{equation*}
\big\| n(t)  \big\|_{H^{1/2}} +\big\|
\partial_t n(t) \big\|_{H^{-1/2}}  \lesssim \exp(c|t|\big\|u_0\big\|^{4m-2}_{L^2}) \max(\big\| n_0  \big\|_{H^{1/2}} +\big\|
n_1 \big\|_{H^{-1/2}}, \big\|u_0\big\|^{2m}_{L^2}).
\end{equation*}
\end{theorem}

The paper is organized as follows.

In Section 2, we first give some scaling analysises of the
criticality, then give the linear and nonlinear estimates along with
 Ginibre, Tsutsumi and Velo  \cite{GiTV97}  in the $X^{s, b}$ spaces, which
 was introduced by Bourgain \cite{Bo93},
 Kenig, Ponce and Vega \cite{KePV93b}, Klainerman and Machedon
\cite{KlM93a}, \cite{KlM95a}. We also can refer to Foschi
\cite{Fo00}, Gr\"{u}nrock \cite{Gr02} and Selberg \cite{Se99} .

For the free dispersive equation of the form
 \begin{equation}\label{equ6}
 iu_t + \varphi (D_x) u = 0, \quad D_x=- i \partial_{x},
 \end{equation}
where $\varphi$ is a measurable function, let $X^{s, b}_{\varphi}$
be the completion of $\mathcal{S}(\mathbb{R})$ with respect to
\begin{equation*}
\aligned
  \|f\|_{X_{\varphi}^{s, b}} : = & \| <\xi>^s <\tau>^b \mathcal{F}(e^{-it \varphi(D_x)} f(x,t))
  \|_{L^2_{\xi, \tau}} \\
   = & \| <\xi>^s <\tau - \varphi(\xi)>^b  \widehat{f}(\xi, \tau))
  \|_{L^2_{\xi, \tau}}
\endaligned
\end{equation*}
In general, we use the notation $X^{s, b}_{\pm}$ for $\varphi(\xi)=
\pm <\xi>$ and $X^{s, b}$ for $\varphi(\xi)=-|\xi|^2$ without
confusion. For a given time interval $I$, we define
\begin{equation*} \aligned \|f\|_{X^{s, b}(I)}
=& \displaystyle\inf_{\widetilde{f}_{|I}= f}
\|\widetilde{f}\|_{X^{s, b}} \quad  \text{where} \quad
\widetilde{f} \in X^{s, b};  \\
\|f\|_{X_{\pm}^{s, b}(I)}= &\displaystyle\inf_{\widetilde{f}_{|I}=
f} \|\widetilde{f}\|_{X^{s, b}_{\pm}} \quad  \text{where} \quad
\widetilde{f} \in X^{s, b}_{\pm}. \endaligned
\end{equation*}

In Section 3, we transform the KLS$_{m}$ $(\ref{equ1})$ into an
equivalent system of first order in $t$ in the usual way, then make
use of Strichartz type estimates to give the local well-posedness in
the $X^{0, b}([0, \delta])\times X^{1/2, b}_{\pm}([0, \delta])$
spaces for some $0< b <1/2$, which is useful for the iteration
procedure. In general, we can obtain the local well-posedness for
$b\geq \frac{1}{2}$,  but in order to get the global wellposedness,
we use $0< b <1/2$ to obtain some gains.

In Section 4, we show that the local result can be iterated to get a
solution on any time interval $[0, T]$. We first can construct the
solution step by step on some time intervals, which is only
dependent of $\big\| u(t) \big\|_{L^2}=\big\| u_0 \big\|_{L^2}$.
Then we can repeat this entire procedure to get the desired time
$T$, each time advancing a time of length $\sim 1/ \big\|
u_0\big\|^{4m-2}_{L^2}$ (independent of $\big\|\big( n(t),
\partial_t n(t)\big)\big\|_{H^{1/2}, H^{-1/2}}$).

We use the following standard facts about the space $X^{s,
b}_{\varphi}$ \cite{Gr02}.

Let $\psi\in C^{\infty}_{0}(\mathbb{R})$ and satisfy $\text{supp} \{
\psi \} \subset (-2, 2)$; $\psi|_{[-1, 1]} = 1$; $\psi(t) =
\psi(-t), \psi \geq 0$. For $0< \lambda \leq 1$, define
$\psi_{\lambda}(t)= \psi(\frac{t}{\lambda})$.

For $s\in \mathbb{R}, b\geq 0$, we have the following homogeneous
estimate
\begin{equation}\label{ine1}
\big\|\psi_{\delta}\ e^{i\varphi(D_x)t} f(x)
\big\|_{X^{s,b}_{\varphi}} \leq c \delta^{\frac{1}{2}-b} \big\| f
\big\|_{H^{s}_{x}};
\end{equation}

\begin{equation}\label{ine6}
\big\| e^{i\varphi(D_x)t} f(x) \big\|_{C(R, H^s_x)} =  \big\| f
\big\|_{H^{s}_{x}};
\end{equation}

For $b'+1 \geq b \geq 0 \geq b' > -\frac{1}{2}$, we have the
following inhomogeneous estimates

\begin{equation}\label{ine3}
\bigg\|\psi_{\delta}\displaystyle \int^{t}_{0}
e^{i(t-s)\varphi(D_x)} F(s)ds \bigg\|_{X^{s,b}_{\varphi}} \leq  c
\delta^{1+b'-b} \big\|F\big\|_{X^{s,b'}_{\varphi}};
\end{equation}

\begin{equation}\label{ine2}
 \bigg\|\displaystyle \int^{t}_{0} e^{i(t-s)\varphi(D_x)}
F(s)ds \bigg\|_{C([0, \delta], H^s_x)} \leq  c
\delta^{\frac{1}{2}+b'} \big\|F\big\|_{X^{s,b'}_{\varphi}};
\end{equation}


For $1< p \leq 2, b\leq \frac{1}{2}- \frac{1}{p}$, we have the
following Sobolev inequality
\begin{equation}\label{ine4}
\big\| f \big\|_{X^{s,b}_{\varphi}} \leq c \big\|
f\big\|_{L^p_t(\mathbb{R}, H^s_x(\mathbb{R}))}.
\end{equation}

Last we introduce the following notation: For $\lambda \in
\mathbb{R}$, Japanese symbol $<\lambda>$ denotes $\big(1+|\lambda|
^2\big)^{1/2}$; $a+$ (resp. $a-$) denotes a number slightly larger
(resp. smaller) than $a$.

\section{Linear and Nonlinear Estimates}
\setcounter{section}{2}\setcounter{equation}{0} In this section, we
first transform the Klein-Gordon-Schr\"{o}dinger system into an
equivalent system of first order in $t$ in the usual way to discuss
the notion of criticality for the system $(\ref{equ1})$. Later, we
give some useful linear and nonlinear estimates.

First, for the notion of criticality, we define
\begin{equation*}
n_{\pm} := \frac{1}{2} \big( n \pm \frac{1}{iA} n_t\big), \quad A
=(I-\partial^2_{x})^{\frac{1}{2}} .
\end{equation*}
Then we have
\begin{equation*}
n = n_{+} + n_{-},\quad   n_t =   i A( n_{+} - n_{-} ), \quad n_{+}
= \overline{n}_{-}.
\end{equation*}
and the equivalent system is
\begin{equation}{\label{equ2}}
\left\{ \aligned
iu_t + \partial^2_x u =& -m(n_{+}+n_{-})|u|^{2(m-1)}u\\
 i \partial_t n_{\pm } \pm A n_{\pm } =& \pm
 \dfrac{1}{2} A^{-1} \big(|u|^{2m} \big)\\
u(0)=u_0 \in H^k_x , &\quad  n_{\pm}(0) = \displaystyle
\frac{1}{2}\big(n_0 \pm
 \frac{1}{iA} n_1\big)\in H^l_x.
\endaligned
\right.
\end{equation}

We follow with Ginibre et at. \cite{GiTV97} to discuss the
criticality through scaling. Consider the following similar system
\begin{equation}\label{equ3}
\left\{ \aligned
i\partial_t u + \partial^2_x u =& -m(n_{+}+n_{-})|u|^{2(m-1)}u\\
 i \partial_t n_{\pm } \pm  (-\partial^2_x)^{1/2}n_{\pm }  =& \pm
 \dfrac{1}{2} (-\partial^2_x)^{-1/2} \big(|u|^{2m} \big).
\endaligned
\right.
\end{equation}

If there were not the term $  \partial^2_x u$ in the LHS of the
first equation in $(\ref{equ3})$, then the system $(\ref{equ3})$
would be invariant under the dilation
\begin{equation*}
\aligned u\rightarrow &u_\lambda=\lambda^{3/(4m-2)}u(\lambda t,
\lambda
x)\\
n\rightarrow & n_\lambda=\lambda^{(2-m)/(2m-1)}n(\lambda t, \lambda
x)
\endaligned
\end{equation*}
and the system $(\ref{equ3})$ would be critical for $(u_0,
n_{\pm}(0))\in H^k_x \times H^l_x$ for
$k=\frac{d}{2}-\frac{3}{4m-2}$, $l=\frac{d}{2}-\frac{2-m}{2m-1}$.
Hence it is $L^2_x\times H^{1/2}_x$-subcritical case for $d=1, 1
\leq m <2$.

If there were not the term $i \partial_t u$ in the LHS of the first
equation in $(\ref{equ3})$, then the system $(\ref{equ3})$ would be
invariant under the dilation
\begin{equation*}
\aligned u\rightarrow &u_\lambda=\lambda^{2/(2m-1)}u(\lambda t,
\lambda
x)\\
n\rightarrow & n_\lambda=\lambda^{2/(2m-1)}n(\lambda t, \lambda x)
\endaligned
\end{equation*}
and the system $(\ref{equ3})$ would be critical for $(u_0,
n_{\pm}(0))\in H^k_x \times H^l_x$ for
$k=\frac{d}{2}-\frac{2}{2m-1}$, $l=\frac{d}{2}-\frac{2}{2m-1}$.
Hence it is $L^2_x\times L^2_x$-subcritical case for $d=1, 1 \leq m
<\frac{5}{2}$.

If there were not the term $\pm  (-\partial^2_x)^{1/2}n_{\pm }$ in
the LHS of the second equation in $(\ref{equ3})$, then the system
$(\ref{equ3})$ would be invariant under the dilation
\begin{equation*}
\aligned u\rightarrow &u_\lambda=\lambda^{5/(4m-2)}u(\lambda^2 t,
\lambda
x)\\
n\rightarrow & n_\lambda=\lambda^{(3-m)/(2m-1)}n(\lambda^2 t,
\lambda x)
\endaligned
\end{equation*}
and the system $(\ref{equ3})$ would be critical for $(u_0,
n_{\pm}(0))\in H^k_x \times H^l_x$ for
$k=\frac{d}{2}-\frac{5}{4m-2}$, $l=\frac{d}{2}-\frac{3-m}{2m-1}$.
Hence it is $L^2_x\times H^{1/2}_x$-subcritical case for $d=1, 1
\leq m <3$.

If there were not the term $i \partial_t n_{\pm }$ in the LHS of the
second equation in $(\ref{equ3})$, then the system $(\ref{equ3})$
would be invariant under the dilation
\begin{equation*}
\aligned u\rightarrow &u_\lambda=\lambda^{2/(2m-1)}u(\lambda^2 t,
\lambda
x)\\
n\rightarrow & n_\lambda=\lambda^{2/(2m-1)}n(\lambda^2 t, \lambda x)
\endaligned
\end{equation*}
and the system $(\ref{equ3})$ would be critical for $(u_0,
n_{\pm}(0))\in H^k_x \times H^l_x$ for
$k=\frac{d}{2}-\frac{2}{2m-1}$, $l=\frac{d}{2}-\frac{2}{2m-1}$.
Hence it is $L^2_x\times L^2_x$-subcritical case for $d=1, 1 \leq m
<\frac{5}{2}$.

That is the reason why we here focus on the local and global
wellposedness of the system $(\ref{equ2})$ in $L^2(\mathbb{R})\times
H^{1/2}(\mathbb{R})$ for $1 \leq m < 2$. We will take the other
cases into account in the forthcoming papers.

Second, we give some known linear estimates. for the Schr\"{o}dinger
equation, we have
\begin{lemma}[Strichartz estimate]\cite{Gr02}\label{lem2-1}
Assume that $ 4 \leq q \leq +\infty, 2 \leq r \leq +\infty$, $0 \leq
\frac{2}{q}\leq \frac{1}{2} -\frac{1}{r} $, and
$$s= \frac{1}{2}
-\frac{1}{r}-\frac{2}{q}.$$ Then we have
\begin{equation*}
 \big\| u \big\|_{L^q_tL^r_x(\mathbb{R})} \leq c \big\| u
\big\|_{X^{s, \frac{1}{2}+}}.
\end{equation*}
\end{lemma}

In particular, combining with the trivial equality $\big\|u
\big\|_{L^2_{t,x}} =\big\| u \big\|_{X^{0,0}}$, we have

\begin{lemma}\cite{Gr02}\label{lem2}
Assume that $0 < \frac{1}{r} \leq \frac{1}{2}$,
$\frac{1}{2}-\frac{1}{r} \leq \frac{2}{q} < \frac{1}{2} +
\frac{1}{r}$ and
\begin{equation*}
\begin{array}{rl}
b> \frac{1}{2}- \frac{1}{q}+ \frac{1}{2}(\frac{1}{2}-\frac{1}{r}).
\end{array}
\end{equation*}
Then the estimate
\begin{equation}
\begin{array}{rl}
\big\|u \big\|_{L^q_tL^r_x(\mathbb{R})}\leq c \big\| u
\big\|_{X^{0,b}}
\end{array}
\end{equation}
holds true for all $u\in X^{0,b}$.
\end{lemma}

For the Klein-Gordon equation, we will use the fact that
\begin{equation}\label{ine10}
\big\| n_{\pm} \big\|_{L^{p}_{t}H^{s}_x(\mathbb{R})} \leq c \big\|
n_{\pm} \big\|_{X^{s, b}},\ \text{for}\  2<p<\infty,
b>\frac{1}{2}-\frac{1}{p}
\end{equation}
which can be obtained from the interpolation between
\begin{equation*}
\big\| n_{\pm} \big\|_{L^{\infty}_{t}H^{s}_x} \leq c \big\| n_{\pm}
\big\|_{X^{s, \frac{1}{2}+}_{\pm}}
\end{equation*}
and the trivial equality
\begin{equation*}
\big\| n_{\pm} \big\|_{L^{2}_{t}H^{s}_x} = \big\| n_{\pm}
\big\|_{X^{s, 0}_{\pm}}.
\end{equation*}

%

Finally, we give some useful nonlinear estimates, which are
especially important to the iteration procedure.
\begin{lemma}[Nonlinear estimates]\label{lem4}
Let $1 \leq m <2$, then there exists some $0< \epsilon <
1-\frac{m}{2}$, such that the estimates
\begin{equation*}
\aligned \big\| n_{\pm} \big|u\big|^{2(m-1)} u\big\|_{X^{0, b'_1}}
&\leq c \big\|n_{\pm} \big\|_{X^{\frac{1}{2}, b_2}_{\pm}} \big\|u
\big\|^{2m-1}_{X^{0, b_1}} \\
\big\| \big|u\big|^{2m} \big\|_{X^{-\frac{1}{2}, b'_2}_{\pm}} &\leq
c \big\| u\big\|^{2m}_{X^{0, b_1}}
\endaligned
\end{equation*}
hold for any $n_{\pm} \in X^{\frac{1}{2}, b_2}_{\pm}$ and $u \in
X^{0, b_1}$ where $b_1=b_2=\frac{2m-1}{4m}+\epsilon$,
$b'_1=b'_2=-\frac{1}{2}+2m \epsilon$.
\end{lemma}

{\bf Proof:} By $(\ref{ine4})$, H\"{o}lder inequality,  we have
\begin{equation*}
\aligned \big\| n_{\pm} \big|u\big|^{2(m-1)} u\big\|_{X^{0, b'_1}}
&\leq c \big\| n_{\pm} \big|u\big|^{2(m-1)}
u\big\|_{L^{\frac{1}{1-2m\epsilon }}_t
L^{2}_x} \\
&\leq c \big\|n_{\pm} \big\|_{L^{4m}_tL^{\frac{1}{\theta}}_x}
\big\|u
\big\|^{2m-1}_{L^{q}_tL^{\frac{4m-2}{1-2\theta}}_x}\\
&\leq c \big\|n_{\pm} \big\|_{L^{4m}_tH^{1/2}_x} \big\|u
\big\|^{2m-1}_{L^{q}_tL^{\frac{4m-2}{1-2\theta}}_x}
\endaligned
\end{equation*}
where $\frac{1}{4m}+\frac{2m-1}{q}=1-2m\epsilon$ and $0 < \theta \ll
1$.

From $(\ref{ine10})$, we have
$$\big\|n_{\pm} \big\|_{L^{4m}_tH^{1/2}_x} \leq c \big\|n_{\pm} \big\|_{X^{\frac{1}{2}, b_2}_{\pm}};$$

From Lemma $\ref{lem2}$, we have
$$\big\|u
\big\|_{L^{q}_tL^{\frac{4m-2}{1-2\theta}}_x} \leq c  \big\|u
\big\|_{X^{0, b_1}},$$ under the conditions
\begin{equation}\label{ine11}
\left\{
\aligned \theta+2\epsilon & <2-m,\\
\theta-4m\epsilon & <m+\frac{1}{2m}-2,
\endaligned
\right.
\end{equation}
which can be satisfied for $1\leq m < 2$.

Therefore, we obtain
\begin{equation*}
\big\| n_{\pm} \big|u\big|^{2(m-1)} u\big\|_{X^{0, b'_1}} \leq c
\big\|n_{\pm} \big\|_{X^{\frac{1}{2}, b_2}_{\pm}} \big\|u
\big\|^{2m+1}_{X^{0, b_1}}. \\
\end{equation*}

In addition, by $(\ref{ine4})$ and Sobolev inequality and Lemma
$\ref{lem2}$, we have
\begin{equation*}
\aligned \big\| \big|u\big|^{2m} \big\|_{X^{-\frac{1}{2},
b'_2}_{\pm}} \leq & c \big\| \big|u\big|^{2m}
\big\|_{L^{\frac{1}{1-2m\epsilon}}_tH^{-1/2}_x}\\
\leq & c \big\| \big|u\big|^{2m}
\big\|_{L^{\frac{1}{1-2m\epsilon}}_tL^{1}_x} \\
\leq & c \big\| u
\big\|^{2m}_{L^{\frac{2m}{1-2m\epsilon}}_tL^{2m}_x} \\
 \leq & c \big\| u\big\|^{2m}_{X^{0, b_1}}.
\endaligned
\end{equation*}

The proof is completed.

\begin{remark}
If $1\leq m \leq 1+\frac{\sqrt{2}}{2}$, we take the value of
$\theta$ and $\epsilon$ in the region $ABC$ of Figure $1$; If
$1+\frac{\sqrt{2}}{2}\leq m \leq 1+ \frac{\sqrt{3}}{2}$, we take
the value of $\theta$ and $\epsilon$ in the region $ABoD$ of
Figure $2$ (see next page); If $1+\frac{\sqrt{3}}{2} \leq m < 2$,
we take the value of $\theta$ and $\epsilon$ in the region $ABo$
of Figure $3$ (see next page). As we know, when $m=2$, it is
difficult to prove the nonlinear estimates in the above lemma for
some $b_1, b_2, b'_1, b'_2$ satisfying
$2m+b'_1+b'_2=(4m-1)b_1+b_2$. Hence we cannot prove the global
well-posedness in Theorem $\ref{the1}$ for the endpoint case
$m=2$.
\end{remark}

\begin{figure}
\centering
\includegraphics[width=0.6\textwidth]{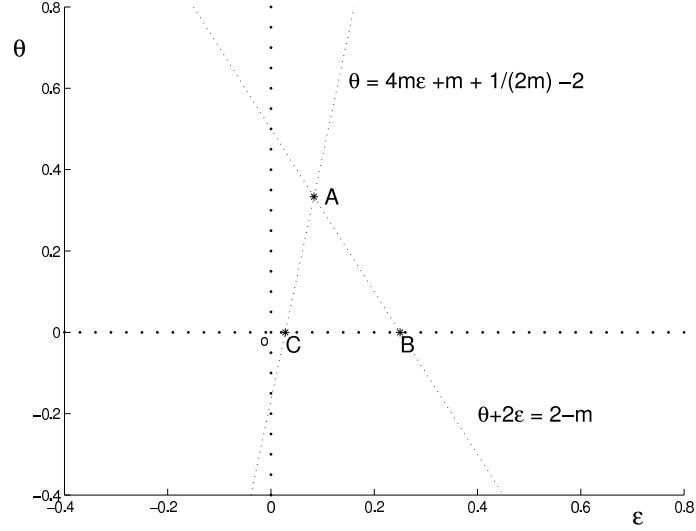}
\caption[]{ $\theta-\epsilon$ parameter picture for $1\leq m \leq
1+\frac{\sqrt{2}}{2}.$}
\end{figure}
\begin{figure}
\centering
\includegraphics[width=0.6\textwidth]{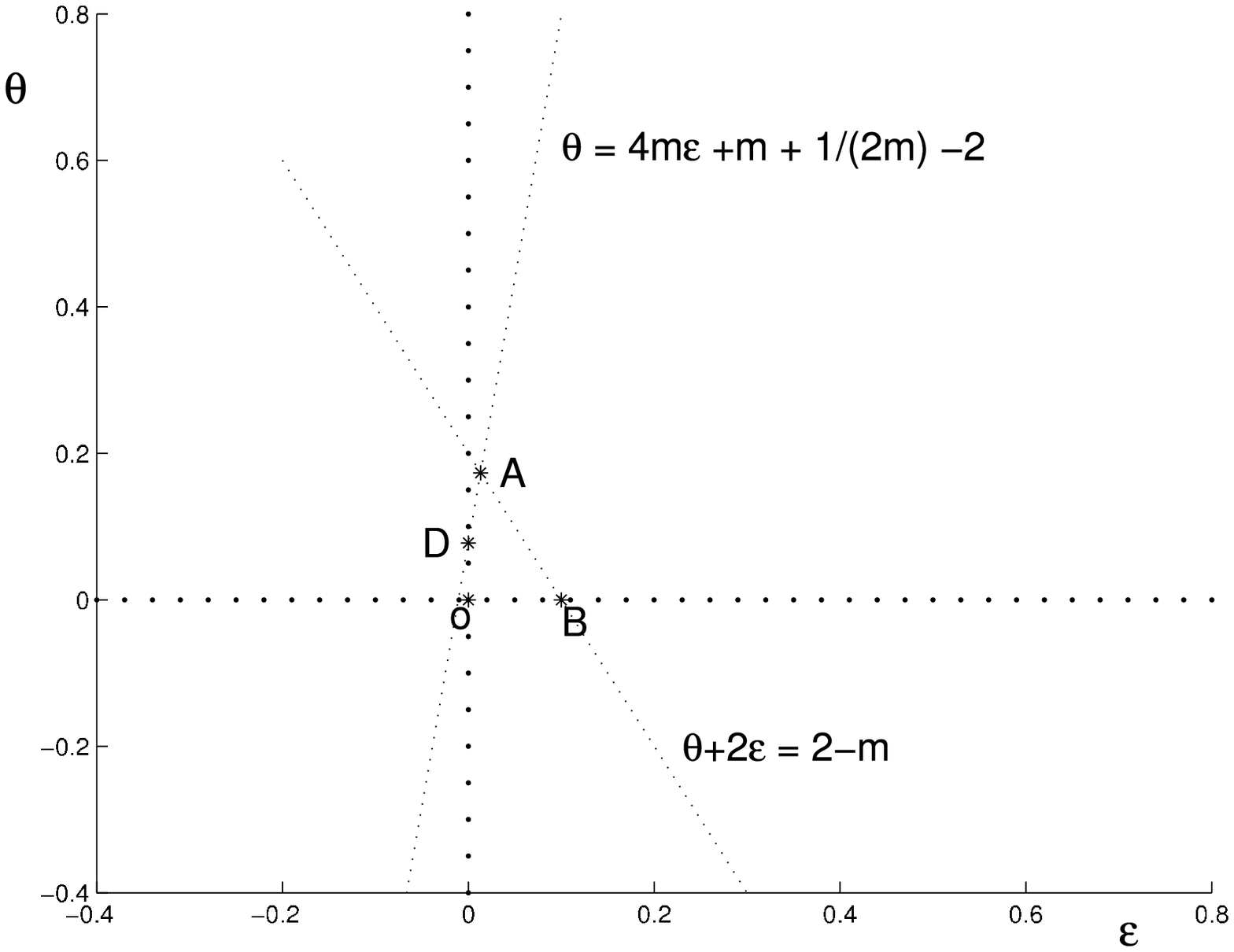}
\caption[]{ $\theta-\epsilon$ parameter picture for
$1+\frac{\sqrt{2}}{2}\leq m \leq 1+\frac{\sqrt{3}}{2}.$}
\end{figure}
\begin{figure}
\centering
\includegraphics[width=0.6\textwidth]{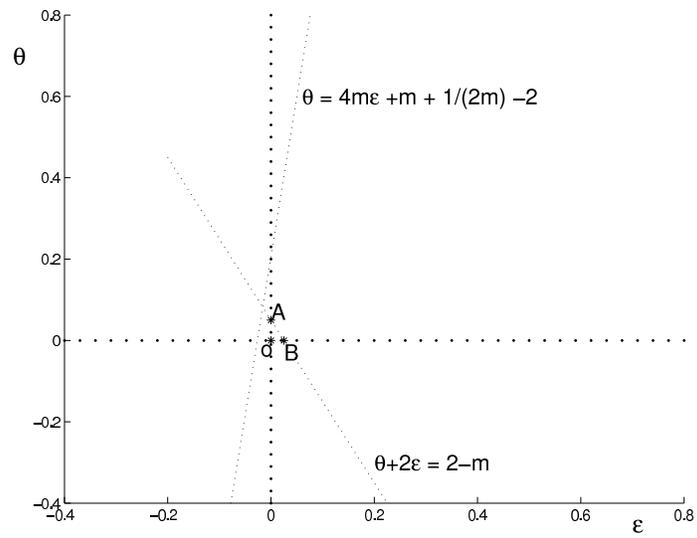}
\caption[]{ $\theta-\epsilon$ parameter picture for
$1+\frac{\sqrt{3}}{2} \leq m < 2.$}
\end{figure}

\section{The Local Well-Posedness}
\setcounter{section}{3} \setcounter{equation}{0}In this section, We
construct a solution of $(\ref{equ2})$ in some time interval $[0,
\delta]$ using the fixed point argument.

The KLS$_m$ $(\ref{equ2})$ has the following equivalent integral
equation formulation
\begin{equation*}
\aligned u(t) &= U(t)u_0 + imU_{*R}[(n_{+}+n_{-})|u|^{2(m-1)}u](t);\\
n_{\pm}&=W_{\pm}n_{\pm}(0) \mp \frac{i}{2}W_{\pm *
R}(A^{-1}|u|^{2m})(t).
\endaligned
\end{equation*}
where
\begin{equation*}
\aligned U(t)u_0 = e^{it\partial^2_x}u_0,\ \ \ \ \ \  \quad & \ \ \ U_{*R}F(t)= \int^t_0 U(t-s)F(s)ds;\\
W_{\pm}(t)n_{\pm}(0)=e^{\pm itA}n_{\pm}(0), \quad & W_{\pm
* R}G(t)= \int^t_0 W(t-s)G(s)ds.
\endaligned
\end{equation*}

For $0<\delta<1$, we define a mapping $M=\big(\Lambda_S(u, n_{\pm}),
\Lambda_{W_{\pm}}(u, n_{\pm})\big)$ by
\begin{equation*}
\left\{
\aligned \Lambda_S(u,
n_{\pm}) &= \psi_{\delta} U(t)u_0 + im\psi_{\delta} U_{*R}[(n_{+}+n_{-})|u|^{2(m-1)}u](t);\\
\Lambda_{W_{\pm}}(u, n_{\pm}) &=\psi_{\delta} W_{\pm}n_{\pm}(0) \mp
\frac{i}{2} \psi_{\delta}  W_{\pm * R}(A^{-1}|u|^{2m})(t).
\endaligned
\right.
\end{equation*}

\begin{proposition}[\bf The local well-posedness]{\label{pro1}}
Let $1 \leq m <2$, $\epsilon>0$ as in Lemma $\ref{lem4}$,
$b_1=b_2=\frac{2m-1}{4m}+\epsilon$, $b'_1=b'_2=-\frac{1}{2}+2m
\epsilon$. Assume that $u_0 \in L^2(\mathbb{R})$, $n_{\pm}(0) \in
H^{1/2}(\mathbb{R})$. Then there exists a positive number $\delta$
satisfying
\begin{eqnarray}
\delta^{m+\frac{1}{2}+b'_2-(2m-1)b_1-b_2} \big\| u_0
\big\|^{2m-1}_{L^2} & \lesssim & 1 ; \label{ine5}\\
\delta^{m+\frac{1}{2}+b'_1-(2m-1)b_1-b_2} \big\| u_0
\big\|^{2m-2}_{L^2} \big\| n_{\pm}(0) \big\|_{H^{1/2}} & \lesssim & 1 ; \label{ine8}\\
\delta^{m+\frac{1}{2}+b'_2-2mb_1} \big\| u_0 \big\|^{2m}_{L^2} &
\lesssim &  \big\| n_{\pm}(0) \big\|_{H^{1/2}}, \label{ine9}
\end{eqnarray}
such that the above Cauchy problem $(\ref{equ2})$ has a unique
solution $u(t,x) \in C([0, \delta], L^2)$ and $n_{\pm}(t,x) \in
C([0, \delta], H^{1/2})$ with the property
\begin{equation*}
\big\|u\big\|_{X^{0, b_1}([0, \delta])} \lesssim
\delta^{\frac{1}{2}-b_1} \big\| u_0\big\|_{L^2}, \quad \big\|
n_{\pm}\big\|_{X^{\frac{1}{2}, b_2}_{\pm}([0, \delta])} \lesssim
\delta^{\frac{1}{2}-b_2} \big\| n_{\pm}(0) \big\|_{H^{1/2}}.
\end{equation*}
\end{proposition}

\begin{remark}
According to the value of $b_1, b_2, b'_1, b'_2$, we have
\begin{equation*}
\aligned
m+\frac{1}{2}+b'_2-(2m-1)b_1-b_2&=m+\frac{1}{2}+b'_1-(2m-1)b_1-b_2\\
&=m+\frac{1}{2}+b'_2-2mb_1\\
&=\frac{1}{2}.
\endaligned
\end{equation*}
\end{remark}

{\bf Proof: } Define the closed set $Y$ as
$$ Y = \bigg\{  \big\|u\big\|_{X^{0, b_1}([0, \delta])} \leq 2 c
\delta^{\frac{1}{2}-b_1} \big\| u_0\big\|_{L^2}, \quad \big\|
n_{\pm}\big\|_{X^{\frac{1}{2}, b_2}_{\pm}([0, \delta])} \leq 2 c
\delta^{\frac{1}{2}-b_2} \big\| n_{\pm}(0) \big\|_{H^{1/2}}.
\bigg\}$$  Define the metric in the set $Y$ as
$$d\big( (u_1, n_{1\pm}), (u_2, n_{2\pm})  \big)=  \big\|  u_1 -u_2\big\|_{X^{0, b_1}([0, \delta])}+ \big\|n_{1\pm} -n_{2\pm}  \big\|_{X^{\frac{1}{2}, b_2}_{\pm}([0, \delta])} .$$

First, we prove that $M$ maps $Y$ into itself under some conditions
on $\delta$.

Now take any $(u, n_{\pm}) \in Y$. By $(\ref{ine1})$, $(\ref{ine3})$
and Lemma $\ref{lem4}$, we have
\begin{equation*}
\aligned \big\|  \Lambda_S(u, n_{\pm})  \big\|_{X^{0, b_1}([0,
\delta])} &\leq c \delta^{\frac{1}{2}-b_1} \big\| u_0\big\|_{L^2} +
c \delta^{1+b'_1-b_1} \big\| n_{\pm} \big|u\big|^{2(m-1)}
u\big\|_{X^{0,
b'_1}} \\
& \leq c \delta^{\frac{1}{2}-b_1} \big\| u_0\big\|_{L^2} + c
\delta^{1+b'_1-b_1} \big\|n_{\pm} \big\|_{X^{\frac{1}{2},
b_2}_{\pm}} \big\|u \big\|^{2m-1}_{X^{0, b_1}} \\
& \leq c \delta^{\frac{1}{2}-b_1} \big\| u_0\big\|_{L^2} + c
\delta^{1+b'_1-b_1} \delta^{\frac{1}{2}-b_2} \big\| n_{\pm}(0)
\big\|_{H^{1/2}} \big( \delta^{\frac{1}{2}-b_1} \big\|
u_0\big\|_{L^2} \big)^{2m-1} \\
& \leq 2 c \delta^{\frac{1}{2}-b_1} \big\| u_0\big\|_{L^2}
\endaligned
\end{equation*}
under the condition  $(\ref{ine8})$.

In addition, we have
\begin{equation*}
\aligned \big\| \Lambda_{W_{\pm}}(u, n_{\pm})
\big\|_{X^{\frac{1}{2}, b_2}_{\pm}([0, \delta])} &\leq c
\delta^{\frac{1}{2}-b_2} \big\| n_{\pm}(0) \big\|_{H^{1/2}} + c
\delta^{1+b'_2-b_2}  \big\||u|^{2m}\big\|
_{X^{-\frac{1}{2}, b'_2}_{\pm}} \\
&\leq c \delta^{\frac{1}{2}-b_2} \big\| n_{\pm}(0) \big\|_{H^{1/2}}
+ c \delta^{1+b'_2-b_2} \big\| u\big\|^{2m}_{X^{0, b_1}}\\
&\leq c \delta^{\frac{1}{2}-b_2} \big\| n_{\pm}(0) \big\|_{H^{1/2}}
+ c \delta^{1+b'_2-b_2} \big( \delta^{\frac{1}{2}-b_1} \big\|
u_0\big\|_{L^2} \big)^{2m} \qquad \qquad \qquad \quad\\
&\leq 2 c \delta^{\frac{1}{2}-b_2} \big\| n_{\pm}(0)
\big\|_{H^{1/2}}
\endaligned
\end{equation*}
under the condition $(\ref{ine9})$. Therefore, we prove that $M$
maps $Y$ into itself.

Second, we can prove that $M$ is a contraction map under another
conditions on $\delta$.

Take any $(u_1, n_{1\pm}), (u_2, n_{2\pm}) \in Y$, we have
\begin{equation*}
\aligned &  \big\|  \Lambda_S(u_1, n_{1\pm}) -  \Lambda_S(u_2,
n_{2\pm}) \big\|_{X^{0, b_1}([0, \delta])} \\ \leq & c
\delta^{1+b'_1-b_1} \big\| n_{1\pm} \big|u_1\big|^{2(m-1)} u_1 -
n_{2\pm} \big|u_2\big|^{2(m-1)} u_2 \big\|_{X^{0,
b'_1}} \\
 \leq  & c \delta^{1+b'_1-b_1} \bigg( \big\|n_{1\pm}
\big\|_{X^{\frac{1}{2}, b_2}_{\pm}} \big( \big\|u_1
\big\|^{2(m-1)}_{X^{0, b_1}} + \big\|u_2 \big\|^{2(m-1)}_{X^{0,
b_1}} \big) \big\|u_1 -u_2
\big\|_{X^{0, b_1}} \\
 &\qquad \qquad \qquad \qquad + \big( \big\|u_1
\big\|^{2m-1}_{X^{0, b_1}} + \big\|u_2 \big\|^{2m-1}_{X^{0, b_1}}
\big)\big\|n_{1\pm} -n_{2\pm} \big\|_{X^{\frac{1}{2},
b_2}_{\pm}} \bigg) \\
 \leq & c \delta^{1+b'_1-b_1} \bigg( \delta^{\frac{1}{2}-b_2} \big\| n_{\pm}(0)\big\|_{H^{1/2}} \delta^{(m-1)(1-2b_1)}\big\|u_0\big\|^{2(m-1)}_{L^2} \big\|u_1 -u_2
\big\|_{X^{0, b_1}} \\
&\qquad \qquad \qquad +
\delta^{(2m-1)(\frac{1}{2}-b_1)}\big\|u_0\big\|^{2m-1}_{L^2}
\big\|n_{1\pm} -n_{2\pm}
\big\|_{X^{\frac{1}{2}, b_2}_{\pm}} \bigg)\\
 \leq & \frac{1}{4} \big(\big\|  u_1 -u_2\big\|_{X^{0, b_1}([0, \delta])}+ \big\|n_{1\pm} -n_{2\pm}  \big\|_{X^{\frac{1}{2}, b_2}_{\pm}([0, \delta])} \big)  \\
\endaligned
\end{equation*}
under the  conditions $(\ref{ine8})$ and
$$ \delta^{m+\frac{1}{2}+b'_1-2mb_1} \big\| u_0
\big\|^{2m-1}_{L^2}  \lesssim  1 $$ which is equivalent to
$(\ref{ine5})$ for $b_1=b_2$ and $b'_1=b'_2$.

In addition, we have
\begin{equation*}
\aligned & \big\| \Lambda_{W_{\pm}}(u_1, n_{1\pm})-
\Lambda_{W_{\pm}}(u_2, n_{2\pm}) \big\|_{X^{\frac{1}{2},
b_2}_{\pm}([0, \delta])}\\
\leq & c \delta^{1+b'_2-b_2} \big\||u_1|^{2m}-|u_2|^{2m}\big\|
_{X^{-\frac{1}{2}, b'_2}_{\pm}} \\
\leq & c \delta^{1+b'_2-b_2} \big(\big\| u_1\big\|^{2m-1}_{X^{0, b_1}} + \big\| u_2\big\|^{2m-1}_{X^{0, b_1}} \big)\big\| u_1 -u_2\big\|_{X^{0, b_1}}\\
\leq & c \delta^{1+b'_2-b_2} \big( \delta^{\frac{1}{2}-b_1} \big\| u_0\big\|_{L^2}\big)^{2m-1}\big\| u_1 -u_2\big\|_{X^{0, b_1}}\\
\leq & \frac{1}{4} \big\| u_1 -u_2\big\|_{X^{0, b_1}}
\endaligned
\end{equation*}
under the condition $(\ref{ine5})$.

The standard fixed point arguments gives a unique solution in time
interval $[0, \delta]$. According to $(\ref{ine6})$ and
$(\ref{ine2})$, we can get that $u \in C([0, \delta], L^2)$ and
$n_{\pm} \in C([0, \delta], H^{1/2})$. Summarizing, The proof is
completed.

\section{Global Well-posedness}
\setcounter{section}{4} \setcounter{equation}{0}In this section, we
show that the process can be iterated to get a solution on any time
interval $[0, T]$. We first can construct the solution step by step
on some time intervals, which is only dependent of $\big\| u(t)
\big\|_{L^2}=\big\| u_0 \big\|_{L^2}$. Then we can repeat this
entire procedure to get the desired time $T$.

According the mass conservation in $(\ref{cl})$, we conclude that
$\big\| u(t)\big\|_{L^2} = \big\| u(0)\big\|_{L^2}$. In order to
iterate the local result to obtain the global well-posedness, we are
only concerned with the growth in $\big\|
n_{\pm}(t)\big\|_{H^{1/2}}$ from one time step to the next step.

Suppose that after some number of iterations we reach a time where
$\big\| n_{\pm}(t) \big\|_{H^{1/2}} \gg \big\| u(t)
\big\|^{2m}_{L^2} =\big\| u_0 \big\|^{2m}_{L^2}$. Take this time
position as the initial time $t=0$ so that $\big\| u_0
\big\|^{2m}_{L^2} \ll  \big\| n_{\pm}(0) \big\|_{H^{1/2}}$. Then
$(\ref{ine9})$ is automatically satisfied and by $(\ref{ine8})$, we
may select a time increment of size
\begin{equation}\label{ts}
\delta\sim \big( \big\| u_0 \big\|^{2m-2}_{L^2} \big\| n_{\pm}(0)
\big\|_{H^{1/2}} \big)^{-1/(m+\frac{1}{2}+b'_1-(2m-1)b_1-b_2)}
\end{equation}

Since
$$n_{\pm}(t)= W_{\pm}n_{\pm}(0) \mp
\frac{i}{2}  W_{\pm * R}(A^{-1}|u|^4)(t),$$ We can apply
$(\ref{ine6})$ and $(\ref{ine3})$ and Proposition $\ref{pro1}$ to
obtain
\begin{equation*}
\aligned \big\| n_{\pm}(\delta)\big\|_{H^{1/2}} \leq & \big\|
n_{\pm}(0)\big\|_{H^{1/2}} + c \delta^{\frac{1}{2}+b'_2}\big\|
\big|u\big|^{2m} \big\|_{X^{-\frac{1}{2}, b'_2}_{\pm}} \\
\leq & \big\| n_{\pm}(0)\big\|_{H^{1/2}} +
c\delta^{\frac{1}{2}+b'_2}
\big\| u\big\|^{2m}_{X^{0, b_1}[0, \delta]}\\
\leq & \big\| n_{\pm}(0)\big\|_{H^{1/2}} + c
\delta^{m+\frac{1}{2}+b'_2-2mb_1} \big\| u_0\big\|^{2m}_{L^2}
\endaligned
\end{equation*}
where $c$ is some fixed constant. From this we can see that we can
carry out $N$ iterations on time intervals each of length
$(\ref{ts})$, where
\begin{equation}\label{ms}
N \sim \frac{\big\|
n_{\pm}(0)\big\|_{H^{1/2}}}{\delta^{m+\frac{1}{2}+b'_2-2mb_1} \big\|
u_0\big\|^{2m}_{L^2}}
\end{equation}
before the quantity $\big\| n_{\pm}(t)\big\|_{H^{1/2}} $ doubles.
The total time we advance after these $N$ iterations, by
$(\ref{ts})$ and $(\ref{ms})$ and $2m+b'_1+b'_2=(4m-1)b_1+b_2$, is
\begin{equation*}
N \delta \sim \frac{\big\|
n_{\pm}(0)\big\|_{H^{1/2}}}{\delta^{m-\frac{1}{2}+b'_2-2mb_1} \big\|
u_0\big\|^{2m}_{L^2}} \sim
\frac{1}{\delta^{2m+b'_1+b'_2-(4m-1)b_1-b_2} \big\|
u_0\big\|^{4m-2}_{L^2}} \sim \frac{1}{ \big\|
u_0\big\|^{4m-2}_{L^2}},
\end{equation*}
which is independent of $\big\|n_{\pm}(t)\big\|_{H^{1/2}}$ .

We can now repeat this entire procedure, each time advancing a time
of length $\sim 1/ \big\| u_0\big\|^{4m-2}_{L^2}$. Upon each
repetition, the size of $\big\|n_{\pm}(t)\big\|_{H^{1/2}}$ will at
most double, giving the exponential-in-time upper bound stated in
Theorem $\ref{the1}$.

This completes the proof of Theorem $\ref{the1}$,

{\bf Acknowledge:} We are deeply grateful to Prof. James Colliander
for his valuable suggestions and discussions.

\begin{center}

\end{center}

\begin{thebibliography}{99}
\addcontentsline{toc}{section}{References}

\bibitem{Ba85}A. Bachelot, Probl\`{e}me de Cauchy pour des syst\`{e}mes
hyperboliques semilin\'{e}aires. Ann. Inst. H. Poincar\`{e}, Analyse
non lin\'{e}aire, 1(1985), 453-478.

\bibitem{Ba89}A. Bachelot, Global existence of large amplitude
solutions for Dirac-Klein-Gordon systems in Minkowski space. Lecture
Notes in Mathematics,  1402(1989), 99-113.

\bibitem{BaC78}J. B. Baillon, J. M. Chadam, The Cauchy problem for
the coupled Schr\"{o}dinger-Klein-Gordon equations, G. M. de La
Penha, L. A. Medeiros (eds.), Contemporary developments in Continuum
Mechanics and Partial Differential Equations, North-Horlland
Publishing Company, 1978, 37-44.

\bibitem{Bi89}P. Biler, Asymptotic behavior of solutions and
universal attractors for a system of nonlinear hyperbolic equations.
Integrable systems and applications. M. Balabane, P. Lochak and C.
Sulem (Eds.) Lecture Notes in Physics. No. 342, 27-30.

\bibitem{Bo93}J. Bourgain, Fourier transform restriction phenomena
for certain lattice subsets and applications to nonlinear evolution
equations I \& II. GAFA, 3(1993), 107-156, 209-262.

\bibitem{Bo98a}J.  Bourgain.  Refinements of Strichartz' inequality and
applications to 2D-NLS with critical nonlinearity.  IMRN,  5(1998),
253-283.

\bibitem{Co98}J. Colliander, Wellposedness for Zakharov systems with
generalized nonlinearity. J. Diff. Equa., 148(1998), 351-363.

\bibitem{CoHT06}J. Colliander, J. Holmer and N. Tzirakis, Low
regularity global well-posedness for the Zakharov and
Klein-Gordon-Schr\"{o}dinger systems. arXiv:math.AP/0603595.

\bibitem{Fo00}D. Foschi. On the regularity of multilinear forms associated to the wave
equation. Dissertation, Princeton Univ., 2000.

\bibitem{FuT75a}I. Fukuda, M. Tsutsumi, On the Yukawa-coupled
Klein-Gordon-Schr\"{o}dinger equations in three space dimensions.
Proc. Japan Acad., 51(1975), 402-405.

\bibitem{FuT75b}I. Fukuda, M. Tsutsumi, On coupled
Klein-Gordon-Schr\"{o}dinger equations, I, Bull. Sci. Engrg. Res.
Lab. Waseda Univ. 69(1975), 51-62.

\bibitem{FuT78}I. Fukuda, M. Tsutsumi, On coupled
Klein-Gordon-Schr\"{o}dinger equations, II,  Jour. Math. Anal. Appl.
66 (1978), 358 - 378.

\bibitem{GiTV97}J. Ginibre, Y. Tsutsumi and G. Velo, On the Cauchy
problem for the Zakharov system. J. Funct. Anal., 151(1997),
384-436.

\bibitem{Gr02} A. Gr\"{u}nrock, New applications of the Fourier
restriction norm method to wellposedness problems for nonlinear
evolution equations. Dissertation Univ. Wuppertal, 2002.

\bibitem{KeT98b}M. Keel and T. Tao, Local and global well-posedness
for wave maps on $\mathbb{R}^{1+1}$ for rough data. IMRN, 21(1998),
1117-1156.

\bibitem{KePV93b}C. Kenig, G. Ponce and L. Vega, The Cauchy problem for the
Kortewegde Vries equation in Sobolev spaces of negative indices.
Duke Math. J., 71(1993), 1-21.

\bibitem{KlM93a}S. Klainerman, M. Machedon, Space-time estimates for
null forms and the local existence theorem. Comm. Pure Appl. Math.,
46(1993), 1221-1268.

\bibitem{KlM95a}S. Klainerman, M. Machedon, Smoothing estimates for
null forms and applications. Duke Math. J., 81(1995), 99-134.

\bibitem{Pe04}H. Pecher, Global solutions of the
Klein-Gordon-Schr\"{o}dinger system with rough data. Differential
and Integral Equations, 17(1-2) (2004), 179-214.

\bibitem{Se99}S. Selberg. Multilinear space-time estimates and  applications to local existence
theory for nonlinear wave equations. Dissertation, Princeton Univ.,
1999.

\bibitem{Tz05b}N. Tzirakis, The Cauchy problem for the
Klein-Gordon-Schr\"{o}dinger system in low dimensions below the
energy space. Commun. Part. Diff. Equa., 30(2005), 605-641.

\end{thebibliography}
\end{document}